\documentclass[11pt]{article}
\usepackage{amsfonts, amsmath, amsthm, verbatim}
\usepackage{amssymb}
\usepackage{graphicx}

\renewcommand{\Re}{\mathfrak{Re\,}}
\renewcommand{\Im}{\mathfrak{Im\,}}

\begin{document}

\title{On the Breiman conjecture}
\author{P\'eter Kevei\thanks{
Center for Mathematical Sciences, Technische Universit\"at M\"unchen,
Boltzmannstra{\ss }e 3, 85748 Garching, Germany,
\texttt{peter.kevei@tum.de}}
\and David M. Mason\thanks{
Department of Applied Economics and Statistics, University of Delaware,
213 Townsend Hall, Newark, DE 19716, USA,
\texttt{davidm@udel.edu}}}
\maketitle

\begin{abstract}
Let $Y_{1},Y_{2},\ldots $ be positive, nondegenerate, i.i.d.~$G$ random
variables, and independently let $X_{1},X_{2},\ldots $ be i.i.d.~$F$ random
variables. In this note we show that whenever $\sum X_{i}Y_{i}/\sum Y_{i}$
converges in distribution to nondegenerate limit for some $F\in \mathcal{F}$%
, in a specified class of distributions $\mathcal{F}$, then $G$ necessarily
belongs to the domain of attraction of a stable law with index less than 1.
The class $\mathcal{F}$ contains those nondegenerate $X$ with a finite
second moment and those $X$ in the domain of attraction of a stable law with
index $1<\alpha <2$.
\end{abstract}

\section{Introduction and results}

Let $Y,Y_{1},\ldots $ be positive, nondegenerate, i.i.d.~random variables
with distribution function [df] $G$, and independently let $X, X_{1},\ldots $ be
i.i.d.~nondegenerate random variables with df $F$. Let $\phi _{X}$ denote
the characteristic function [cf] of $X$. We shall use the notation $Y\in D(\beta)$
to mean that $Y$ is in the domain of attraction of a stable law of index $0<
\beta <1$, and $Y\in D(0)$ will denote that $1-G$ is slowly varying at
infinity. Furthermore $\mathcal{RV}_{\infty}(\rho)$ will signify the class
of positive measurable functions regularly varying at infinity with index $%
\rho$, and $\mathcal{RV}_{0}(\rho)$ the class of positive measurable
functions regularly varying at zero with index $\rho$. In particular, using
this notation $Y\in D(\beta)$, with $0\leq \beta <1$, if and only if $1-G\in 
\mathcal{RV}_{\infty }(-\beta )$. \smallskip

For each integer $n\geq 1$ set 
\begin{equation}
T_{n}=\sum_{i=1}^{n}X_{i}Y_{i}/\sum_{i=1}^{n}Y_{i}.  \label{eq:Tdef}
\end{equation}%
Notice that $\mathbb{E}|X|<\infty $ implies that $T_{n}$ is stochastically
bounded. Theorem 4 of Breiman \cite{Brei} says that $T_{n}$ converges in
distribution along the full sequence $\left\{ n\right\} $ for \textit{every} 
$X$ with finite expectation, and with at least one limit law being
nondegenerate if and only if 
\begin{equation}
Y\in D(\beta) ,\text{ with }0\leq \beta <1.  \label{DB}
\end{equation}%
Let $\mathcal{X}$ denote the class of nondegenerate random variables $X$
with $\mathbb{E}|X|<\infty $ and let $\mathcal{X}_{0}$ denote those $X\in $ $%
\mathcal{X}$ such that $\mathbb{E}X=0.$ At the end of his paper Breiman
conjectured that if for \textit{some} $X\in \mathcal{X}$, $T_{n}$ converges
in distribution to some nondegenerate random variable $T$, written 
\begin{equation}
T_{n}\rightarrow_{d} T,\text{ as } n\rightarrow \infty, \text{ with } T 
\text{ nondegenerate,}  \label{dist}
\end{equation}%
then (\ref{DB}) holds$.$ By Proposition 2 and Theorem 3 of \cite{Brei}, for
any $X\in \mathcal{X}$, (\ref{DB}) implies (\ref{dist}), in which case $T$
has a distribution related to the arcsine law. Using this fact, we see that
his conjecture can restated to be: for any $X\in \mathcal{X}$, (\ref{DB}) is
equivalent to (\ref{dist}). \smallskip

It has proved to be surprisingly challenging to resolve.
Mason and Zinn \cite{MZ} partially verified Breiman's
conjecture. They established that whenever $X$ is nondegenerate and
satisfies $\mathbb{E}|X|^{p}<\infty $ for some $p>2,$ then (\ref{DB}) is
equivalent to (\ref{dist}). In this note we further extend this result.

\medskip \noindent \textbf{Theorem}
\textit{Assume that for some $X \in 
\mathcal{X}_{0}$, $1<\alpha \leq 2$, positive slowly varying function $L$ at
zero and $c>0$, 
\begin{equation}
\frac{-\log \left( \Re \phi _{X}(t)\right) }{\left\vert t\right\vert
^{\alpha }L\left( \left\vert t\right\vert \right) }\rightarrow c, \text{ as }%
t \to 0,  \label{alpha}
\end{equation}%
(in the case $\alpha =2$ we assume that $\liminf_{t\searrow 0}L(t)>0$).
Whenever (\ref{dist}) holds then $Y \in D(\beta)$ for some $\beta \in [0,1)$.
} \medskip

Let $\mathcal{F}$ denote the class of random variables that satisfy the
conditions of the theorem. Applying our theorem in combination with
Proposition 2 and Theorem 3 of \cite{Brei} we get the following corollary.

\medskip \noindent \textbf{Corollary} \textit{Whenever $X-\mathbb{E}X \in 
\mathcal{F}$, (\ref{DB}) is equivalent to (\ref{dist})}. \medskip

\noindent \textbf{Remark 1} It can be inferred from Theorem 8.1.10 of
Bingham et~al.~\cite{BGT} that for $X\in \mathcal{X}_{0}$, (\ref{alpha})
holds for some $1<\alpha <2$, positive slowly varying function $L$ at zero
and $c>0$ if and only if $X$ satisfies $\mathbb{P}\left\{ \left\vert
X\right\vert >x\right\} \sim L( 1/x) x^{-\alpha } c\Gamma(\alpha) \frac{2}{%
\pi } \sin \left( \frac{\pi \alpha }{2}\right)$. Note that a random variable 
$X\in \mathcal{X}_{0}$ in the domain of attraction of a stable law of index $%
1<\alpha <2$ satisfies (\ref{alpha}). Also a random variable $X\in \mathcal{X%
}_{0}$ with variance $0<\sigma^2 <\infty$ fulfills (\ref{alpha}) with $%
\alpha =2$, $L=1$ and $c=\sigma^{2}/2$. \medskip

\noindent \textbf{Remark 2} Consult Kevei and Mason \cite{KM} for a fairly
exhaustive study of the asymptotic distributions of $T_{n}$ along
subsequences, along with relevations of their unexpected properties. \medskip

The theorem follows from the two propositions below. First we need more
notation. For any $\alpha \in (1,2]$ define for $n\geq 1$ 
\begin{equation}
S_{n}(\alpha )=\frac{\sum_{i=1}^{n}Y_{i}^{\alpha }}{\left(
\sum_{i=1}^{n}Y_{i}\right) ^{\alpha }}.  \label{eq:snadef}
\end{equation}

\noindent \textbf{Proposition 1} \textit{Assume that the assumptions of the
theorem hold. Then for some $0<\gamma \leq 1$ 
\begin{equation}
\mathbb{E}S_{n}(\alpha )\rightarrow \gamma, \text{ as }n\rightarrow \infty.
\label{gam}
\end{equation}
} \medskip

The next proposition is interesting in its own right. It is an extension of
Theorem 5.3 by Fuchs et~al.~\cite{FJT}, 
where $\alpha =2$ (see also Proposition 3 of \cite{MZ}).\medskip 

\noindent \textbf{Proposition 2} \textit{If (\ref{gam}) holds with some $%
\gamma \in (0,1]$ then $Y \in D(\beta)$, for some $\beta \in [0,1)$, where $%
-\beta \in (-1,0]$ is the unique solution of 
\begin{equation*}
\mathrm{Beta}(\alpha -1,\beta +1)= \frac{\Gamma (\alpha -1)\Gamma (1+\beta )%
}{\Gamma (\alpha -\beta )}=\frac{1}{\gamma (\alpha -1)}.
\end{equation*}%
In particular, $Y \in D(0)$ for $\gamma =1$.} \smallskip

\textit{Conversely, if $G\in D(\beta)$, $0\leq \beta <1$, then (\ref{gam})
holds with 
\begin{equation*}
\gamma =\frac{\Gamma (\alpha -\beta )}{\Gamma (\alpha )\Gamma (1+\beta )}= 
\frac{1}{(\alpha -1) \mathrm{Beta}(\alpha -1,\beta +1)}.
\end{equation*}%
}

\section{Proofs}

Set for each $n\geq 1$, $R_{i}=Y_{i}/\sum_{l=1}^{n}Y_{l},$ for $i=1,\dots ,n$%
. For notational ease we drop the dependence of $R_{i}$ on $n\geq 1$.
Consider the sequence of strictly decreasing continuous functions $\left\{
\varphi _{n}\right\} _{n\geq 1}$ on $[1,\infty)$ defined by $\varphi _{n}
(y) =\mathbb{E}\left( \sum_{i=1}^{n}R_{i}^{y}\right)$, $y\in [ 1,\infty)$.
Note that each function $\varphi _{n}$ satisfies $\varphi _{n}(1) =1$. By a
diagonal selection procedure for each subsequence of $\left\{ n\right\}
_{n\geq 1}$ there is a further subsequence $\{ n_{k} \} _{k\geq 1}$ and a
right continuous nonincreasing function $\psi$ such that $\varphi _{n_{k}}$
converges to $\psi$ at each continuity point of $\psi$.\medskip

\noindent\textbf{Lemma 1} \textit{Each such function $\psi$ is continuous on $%
(1,\infty)$.}\smallskip

\noindent \textit{Proof} Choose any subsequence $\left\{
n_{k}\right\}_{k\geq1}$ and a right continuous nonincreasing function $\psi$
such that $\varphi_{n_{k}}$ converges to $\psi$ at each continuity point of $%
\psi$ in $(1,\infty)$. Select any $x>1$ and continuity points $x_1,x_2 \in
(1, \infty)$ of $\psi$ such that $1< x_1< x < x_2 <\infty$. Set $\rho_{1}=
x_1-1$ and $\rho_2=x_2-1$. Since $\rho_{2}/\rho_{1}>1$ we get by H\"{o}%
lder's inequality 
\begin{equation*}
\sum_{i=1}^{n_k}R_{i}^{x_{1}}=\sum_{i=1}^{n_k}
R_{i}^{\rho_{1}}R_{i}\leq\left( \sum_{i=1}^{n_k}R_{i}^{\rho_{2}}R_{i}\right)
^{\rho _{1}/\rho_{2}}=\left( \sum_{i=1}^{n_k}R_{i}^{x_{2}}\right) ^{\rho
_{1}/\rho_{2}}.
\end{equation*}
Thus by taking expectations and using Jensen's inequality we get $%
\varphi_{n_k}( x_{1} ) \leq \left( \varphi_{n_k} (x_{2})
\right)^{\rho_{1}/\rho_{2}}.$ Letting $n_{k}\rightarrow\infty$, we have $%
\psi(x_{1} ) \leq \left( \psi (x_{2} ) \right)^{\rho_{1}/\rho_{2}}.$ Since $%
x_{1}< x$ and $x_{2}> x$ can be chosen arbitrarily close to $x$ we conclude
by right continuity of $\psi$ at $x$ that $\psi(x-) =\psi (x+) =\psi(x)$. %
\mbox{$\Box$}\medskip

\noindent \textit{Proof of Proposition 1} For a complex $z$, we use the
notation for the principal branch of the logarithm, $Log\left( z\right)
=\log \left\vert z\right\vert +\imath \arg z$, where $-\pi <\arg z\leq \pi $%
, i.e. $z=\left\vert z\right\vert \exp \left( \imath \arg z\right) .$ We see
that for all $t$ 
\begin{equation*}
\begin{split}
\mathbb{E}\exp \left( \imath tT_{n}\right) & =\mathbb{E}\left(
\prod\limits_{j=1}^{n}\phi _{X}\left( tR_{j}\right) \right)  \\
& =\mathbb{E}\left( \prod\limits_{j=1}^{n}\exp \left( Log\phi _{X}\left(
tR_{j}\right) \right) \right) .
\end{split}%
\end{equation*}%
Since $\mathbb{E}X=0$ we have $\Re \phi _{X}(u)=1-o_{+}(u)$, where $%
o_{+}(u)\geq 0$, and $o_{+}(u)$ and $o_{+}(u)/u\rightarrow 0$ as $%
u\rightarrow 0$; and $\Im \phi _{X}(u)=o(u)$. This when combined with 
\begin{equation*}
\left( \arctan \theta \right) ^{\prime }=\frac{1}{1+\theta ^{2}}
\end{equation*}%
gives as $u\rightarrow 0$, 
\begin{equation*}
\arg \phi _{X}(u)=\arctan \left( \frac{\Im \phi _{X}(u)}{\Re \phi _{X}(u)}%
\right) =o\left( u\right) .
\end{equation*}%
Note that for all $\left\vert u\right\vert >0$ sufficiently small so that $%
\Re \phi _{X}(u)>0$ 
\begin{equation*}
Log\phi _{X}(u)=Log(\Re \phi _{X}(u)+\imath \Im \phi _{X}(u))=\log \Re \phi
_{X}(u)+Log\left( 1+\imath \frac{\Im \phi _{X}(u)}{\Re \phi _{X}(u)}\right) ,
\end{equation*}%
where for the second term 
\begin{equation*}
\Re Log\left( 1+\imath \frac{\Im \phi _{X}(u)}{\Re \phi _{X}(u)}\right) =%
\frac{1}{2}\left( \frac{\Im \phi _{X}(u)}{\Re \phi _{X}(u)}\right)
^{2}\left( 1+o\left( u\right) \right) ,\text{ as }u\rightarrow 0.
\end{equation*}%
Thus for every $\varepsilon >0$ for all $\left\vert t\right\vert >0$
sufficiently small and independent of $n\geq 1$ and $R_{1},\dots ,R_{n}$ 
\begin{equation*}
1-\varepsilon ^{2}t^{2}\leq \cos (\varepsilon t)\leq \Re \left( \exp \left\{
\sum_{j=1}^{n}Log\left( 1+\imath \frac{\Im \phi _{X}(tR_{j})}{\Re \phi
_{X}(tR_{j})}\right) \right\} \right) \leq e^{2^{-1}\varepsilon t^{2}}\leq
1+\varepsilon t^{2}.
\end{equation*}%
Thus we obtain 
\begin{equation*}
\begin{split}
\mathbb{E}\exp \Big\{\sum_{j=1}^{n}\log \Re \phi _{X}(tR_{j})\Big\}\left(
1-\varepsilon ^{2}t^{2}\right) & \leq \mathbb{E}\left( \Re \exp \left(
\imath tT_{n}\right) \right)  \\
& =\Re \mathbb{E}\exp \left( \imath tT_{n}\right)  \\
& \leq \mathbb{E}\exp \Big\{\sum_{j=1}^{n}\log \Re \phi _{X}(tR_{j})\Big\}%
(1+\varepsilon t^{2}).
\end{split}%
\end{equation*}%
We shall show (\ref{alpha}) implies that (\ref{gam}) holds for some $%
0<\gamma \leq 1.$ Now using (\ref{alpha}) we get for any $0<\delta <c$ and
all $\left\vert t\right\vert $ small enough independent of $n\geq 1$, 
\begin{equation*}
\begin{split}
& -\varepsilon t^{2}+\log \mathbb{E}\exp \left( -\left( c+\delta \right)
\left\vert t\right\vert ^{\alpha }\left( \sum_{i=1}^{n}R_{i}^{\alpha
}L\left( \left\vert t\right\vert R_{i}\right) \right) \right) \leq \log
\left( \Re \mathbb{E}\exp \left( \imath tT_{n}\right) \right)  \\
& \leq \varepsilon t^{2}+\log \mathbb{E}\exp \left( -\left( c-\delta \right)
\left\vert t\right\vert ^{\alpha }\left( \sum_{i=1}^{n}R_{i}^{\alpha
}L\left( \left\vert t\right\vert R_{i}\right) \right) \right) .
\end{split}%
\end{equation*}%
Next since $\log s/(1-s)\rightarrow -1$ as $s\nearrow 1$, for all $%
\left\vert t\right\vert $ small enough independent of $n\geq 1$ and $%
R_{1},\dots ,R_{n}$, (keeping mind that $\sum_{i=1}^{n}R_{i}=1$ and $%
1<\alpha \leq 2$) 
\begin{equation*}
\begin{split}
& \log \mathbb{E}\exp \left( -\left( c+\delta \right) \left\vert
t\right\vert ^{\alpha }\left( \sum_{i=1}^{n}R_{i}^{\alpha }L\left(
\left\vert t\right\vert R_{i}\right) \right) \right)  \\
& \geq -\left( 1+\frac{\delta }{2}\right) \mathbb{E}\left( 1-\exp \left(
-\left( c+\delta \right) \left\vert t\right\vert ^{\alpha }\left(
\sum_{i=1}^{n}R_{i}^{\alpha }L\left( \left\vert t\right\vert R_{i}\right)
\right) \right) \right) 
\end{split}%
\end{equation*}%
and 
\begin{equation*}
\begin{split}
& \log \mathbb{E}\exp \left( -\left( c-\delta \right) \left\vert
t\right\vert ^{\alpha }\left( \sum_{i=1}^{n}R_{i}^{\alpha }L\left(
\left\vert t\right\vert R_{i}\right) \right) \right)  \\
& \leq -\left( 1-\frac{\delta }{2}\right) \mathbb{E}\left( 1-\exp \left(
-\left( c-\delta \right) \left\vert t\right\vert ^{\alpha }\left(
\sum_{i=1}^{n}R_{i}^{\alpha }L\left( \left\vert t\right\vert R_{i}\right)
\right) \right) \right) .
\end{split}%
\end{equation*}%
Further since $\left( 1-\exp \left( -y\right) \right) /y\rightarrow 1$ as $%
y\searrow 0$, for all $\left\vert t\right\vert $ small enough independent of 
$n\geq 1$, 
\begin{equation*}
\begin{split}
& -\left( 1+\frac{\delta }{2}\right) \mathbb{E}\left( 1-\exp \left( -\left(
c+\delta \right) \left\vert t\right\vert ^{\alpha }\left(
\sum_{i=1}^{n}R_{i}^{\alpha }L\left( \left\vert t\right\vert R_{i}\right)
\right) \right) \right)  \\
& \geq -\left( 1+\delta \right) \left( c+\delta \right) \left\vert
t\right\vert ^{\alpha }\mathbb{E}\left( \sum_{i=1}^{n}R_{i}^{\alpha }L\left(
\left\vert t\right\vert R_{i}\right) \right) 
\end{split}%
\end{equation*}%
and 
\begin{equation*}
\begin{split}
& -\left( 1-\frac{\delta }{2}\right) \mathbb{E}\left( 1-\exp \left( -\left(
c-\delta \right) \left\vert t\right\vert ^{\alpha }\left(
\sum_{i=1}^{n}R_{i}^{\alpha }L\left( \left\vert t\right\vert R_{i}\right)
\right) \right) \right)  \\
& \leq -\left( 1-\delta \right) \left( c-\delta \right) \left\vert
t\right\vert ^{\alpha }\mathbb{E}\left( \sum_{i=1}^{n}R_{i}^{\alpha }L\left(
\left\vert t\right\vert R_{i}\right) \right) .
\end{split}%
\end{equation*}%
Therefore for all $\left\vert t\right\vert $ small enough independent of $n$%
, 
\begin{equation*}
\begin{split}
& -\varepsilon t^{2}-\left( 1+\delta \right) \left( c+\delta \right)
\left\vert t\right\vert ^{\alpha }\mathbb{E}\left(
\sum_{i=1}^{n}R_{i}^{\alpha }L\left( \left\vert t\right\vert R_{i}\right)
\right)  \\
& \leq \log \left( \Re \mathbb{E}\exp \left( \imath tT_{n}\right) \right)  \\
& \leq \varepsilon t^{2}-\left( 1-\delta \right) \left( c-\delta \right)
\left\vert t\right\vert ^{\alpha }\mathbb{E}\left(
\sum_{i=1}^{n}R_{i}^{\alpha }L\left( \left\vert t\right\vert R_{i}\right)
\right) .
\end{split}%
\end{equation*}%
By the Potter's bound, Theorem 1.5.6 (i) in \cite{BGT}, for all $A>1$ and $%
1<\alpha _{1}<\alpha <\alpha _{2}$, for all $t>0$ small enough independent
of $n\geq 1$,%
\begin{equation}
A^{-1}\sum_{i=1}^{n}R_{i}^{^{\alpha _{2}}}\leq \sum_{i=1}^{n}R_{i}^{\alpha
}L\left( \left\vert t\right\vert R_{i}\right) /L\left( \left\vert
t\right\vert \right) \leq A\sum_{i=1}^{n}R_{i}^{\alpha _{1}}.  \label{bb}
\end{equation}%
We see now that for all $n\geq 1$ and $0<4\varepsilon <c$, appropriate $%
1<\alpha _{1}<\alpha <\alpha _{2}$ and all $\left\vert t\right\vert $ small
enough independent of $n$, 
\begin{equation*}
\begin{split}
& -\varepsilon t^{2}-\left( 1+\varepsilon \right) \left( c+2\varepsilon
\right) \left\vert t\right\vert ^{\alpha }L\left( \left\vert t\right\vert
\right) \mathbb{E}S_{n}\left( \alpha _{2}\right)  \\
& =-\varepsilon t^{2}-\left( 1+\varepsilon \right) \left( c+2\varepsilon
\right) \left\vert t\right\vert ^{\alpha }L\left( \left\vert t\right\vert
\right) \mathbb{E}\left( \sum_{i=1}^{n}R_{i}^{\alpha _{2}}\right)  \\
& \leq \log \left( \Re \mathbb{E}\exp \left( \imath tT_{n}\right) \right)  \\
& \leq \varepsilon t^{2}-\left( 1-\varepsilon \right) \left( c-2\varepsilon
\right) \left\vert t\right\vert ^{\alpha }L\left( \left\vert t\right\vert
\right) \mathbb{E}\left( \sum_{i=1}^{n}R_{i}^{\alpha _{1}}\right)  \\
& =\varepsilon t^{2}-\left( 1-\varepsilon \right) \left( c-2\varepsilon
\right) \left\vert t\right\vert ^{\alpha }L\left( \left\vert t\right\vert
\right) \mathbb{E}S_{n}\left( \alpha _{1}\right) .
\end{split}%
\end{equation*}%
Choose any subsequence $\left\{ n_{k}\right\} _{k\geq 1}$ and a right
continuous nonincreasing function $\psi $ such that $\varphi _{n_{k}}$
converges to $\psi $ at each continuity point of $\psi $, which by Lemma 1 above
is all $\left( 1,\infty \right) $. We see that $\mathbb{E}S_{n_{k}}(\alpha
)\rightarrow \psi (\alpha )$, $\mathbb{E}S_{n_{k}}(\alpha _{1})\rightarrow
\psi (\alpha _{1})$ and $\mathbb{E}S_{n_{k}}(\alpha _{2})\rightarrow \psi
(\alpha _{2})$, where necessarily $0<\psi (\alpha _{2})\leq \psi (\alpha
)\leq \psi (\alpha _{1})\leq 1$. (The case $\psi (\alpha _{1})=0$ cannot
happen, since this would imply that $T$ is degenerate.) We see that for all $%
\left\vert t\right\vert $ sufficiently small independent of $n_{k}\geq 1$, 
\begin{equation*}
-\varepsilon -\left( 1+\varepsilon \right) \left( c+3\varepsilon \right)
\psi \left( \alpha _{2}\right) \leq \log \left( \Re \mathbb{E}\exp \left(
\imath tT_{n_{k}}\right) \right) /\left( \left\vert t\right\vert ^{\alpha
}L\left( \left\vert t\right\vert \right) \right) \leq \varepsilon -\left(
1-\varepsilon \right) \left( c-3\varepsilon \right) \psi \left( \alpha
_{1}\right) ,
\end{equation*}%
where for $\alpha =2$ we use the assumption that in this case $%
\liminf_{t\searrow 0}L\left( t\right) >0$. Since $0<4\varepsilon <c$ can be
made arbitrarily small and $0\leq \psi \left( \alpha _{1}\right) -\psi
\left( \alpha _{2}\right) $ can be made as close to zero as desired, by
letting $n_{k}\rightarrow \infty $, we get that for all $\left\vert
t\right\vert $ sufficiently small 
\begin{equation*}
-\varepsilon -\left( 1+\varepsilon \right) \left( c+4\varepsilon \right)
\psi \left( \alpha \right) \leq \log \left( \Re \mathbb{E}\exp \left( \imath
tT\right) \right) /\left( \left\vert t\right\vert ^{\alpha }L\left(
\left\vert t\right\vert \right) \right) \leq \varepsilon -\left(
1-\varepsilon \right) \left( c-4\varepsilon \right) \psi \left( \alpha
\right) ,
\end{equation*}%
which can happen only if $\psi \left( \alpha \right) $ does not depend on $%
\left\{ n_{k}\right\} $. Thus (\ref{gam}) holds for some $0<\gamma \leq 1$,
namely $\gamma =\psi (\alpha )$. \hfill \mbox{$\Box$}\medskip 

\noindent \textit{Proof of Proposition 2} To begin with, we note that
whenever (\ref{gam}) holds, necessarily $\mathbb{E}Y=\infty $. To see this,
write $D_{n}^{\left( 1\right) }=\max_{1\leq i\leq n}Y_{i}/\left(
\sum_{i=1}^{n}Y_{i}\right) $ and observe that%
\begin{equation*}
\begin{split}
\left( D_{n}^{\left( 1\right) }\right) ^{\alpha } & =\max_{1\leq i\leq n}%
\frac{Y_{i}^{\alpha }}{\left( \sum_{i=1}^{n}Y_{i}\right) ^{\alpha }}\leq
S_{n}(\alpha ) \\
& \leq \max_{1\leq i\leq n}\frac{Y_{i}^{\alpha -1}}{\left(
\sum_{i=1}^{n}Y_{i}\right) ^{\alpha -1}}=\left( D_{n}^{\left( 1\right)
}\right) ^{\alpha -1}.
\end{split}%
\end{equation*}%
From these inequalities it is easy to prove that $\mathbb{E}S_{n}(\alpha
)\rightarrow 0,$ $n\rightarrow \infty ,$ if and only if 
\begin{equation}
D_{n}^{\left( 1\right) }\rightarrow _{P}0,\text{ }n\rightarrow \infty .
\label{zero}
\end{equation}%
Proposition 1 of Breiman \cite{Brei} says that (\ref{zero}) holds if and
only there exists a sequence of positive constants $B_{n}$ converging to
infinity such that%
\begin{equation}
\sum_{i=1}^{n}Y_{i}/B_{n}\rightarrow _{P}1,\text{ }n\rightarrow \infty .
\label{LLN}
\end{equation}%
Since $\mathbb{E}Y<\infty $ obviously implies (\ref{LLN}), it readily
follows that $\mathbb{E}S_{n}(\alpha )\rightarrow 0,$ $n\rightarrow \infty ,$
and thus (\ref{gam}) cannot hold.

We shall first prove the first part of Proposition 2. Following similar
steps as in \cite{MZ} we have that 
\begin{equation*}
\begin{split}
\mathbb{E}\frac{\sum_{i=1}^{n}Y_{i}^{\alpha }}{\left(
\sum_{i=1}^{n}Y_{i}\right) ^{\alpha }}& =n\mathbb{E}\frac{Y_{1}^{\alpha }}{%
\left( \sum_{i=1}^{n}Y_{i}\right) ^{\alpha }} \\
& =\frac{n}{\Gamma (\alpha )}\mathbb{E}\int_{0}^{\infty }Y_{1}^{\alpha
}e^{-t\sum_{i=1}^{n}Y_{i}}t^{\alpha -1}\mathrm{d}t \\
& =\frac{n}{\Gamma (\alpha )}\int_{0}^{\infty }t^{\alpha -1}\mathbb{E}\left(
e^{-tY_{1}}Y_{1}^{\alpha }\right) (\mathbb{E}e^{-tY_{1}})^{n-1}\mathrm{d}t \\
& =:\frac{n}{\Gamma (\alpha )}\int_{0}^{\infty }t^{\alpha -1}\phi _{\alpha
}(t)\phi _{0}(t)^{n-1}\mathrm{d}t.
\end{split}%
\end{equation*}%
Next, assuming (\ref{gam}) and arguing as in the proof of Theorem 3 in \cite%
{Brei} we get 
\begin{equation}
s\int_{0}^{\infty }t^{\alpha -1}\phi _{\alpha }(t)e^{s\log \phi _{0}(t)}%
\mathrm{d}t \rightarrow \gamma \Gamma (\alpha ),\quad s\rightarrow \infty .
\label{eq:sconv}
\end{equation}%
For $y\geq 0,$ let $q(y)$ denote the inverse of $-\log \varphi_{0}(t)$.
Changing the variables to $y=-\log \phi _{0}(t)$ and $t=q(y)$, we get from (%
\ref{eq:sconv}) that 
\begin{equation*}
s\int_{0}^{\infty }\left( q(y) \right) ^{\alpha -1}\phi _{\alpha} \left(
q(y)\right) \exp \left( -sy\right) \mathrm{d}q(y) \rightarrow \gamma \Gamma
(\alpha ),\text{ as } s \rightarrow \infty .
\end{equation*}%
By Karamata's Tauberian theorem, see Theorem 1.7.1$^{\prime }$ on page 38 of 
\cite{BGT}, we conclude that 
\begin{equation*}
v^{-1}\int_{0}^{v}\left( q(x) \right) ^{\alpha -1}\phi_{\alpha} \left( q(x)
\right) \mathrm{d}q(x) \rightarrow \gamma \Gamma (\alpha ),\text{ as }%
v\searrow 0,
\end{equation*}%
which, in turn, by the change of variable $y=q(x)$ gives 
\begin{equation*}
\frac{\int_{0}^{t}y^{\alpha -1}\phi _{\alpha }(y)\mathrm{d}y}{-\log \phi
_{0}(t)}\rightarrow \gamma \Gamma (\alpha ),\text{ as } t\searrow 0.
\end{equation*}%
Now using that $-\log \phi _{0}(t)\sim 1-\phi _{0}(t)$ as $t\rightarrow 0$,
we end up with 
\begin{equation*}
\lim_{t\rightarrow 0}\frac{\int_{0}^{t}y^{\alpha -1}\phi _{\alpha }(y)%
\mathrm{d}y}{1-\phi _{0}(t)}=\gamma \Gamma (\alpha ).
\end{equation*}%
Since $\phi _{\alpha }(y)=\int_{0}^{\infty }e^{-uy}u^{\alpha }G(\mathrm{d}u)$%
, by Fubini's theorem 
\begin{equation*}
\begin{split}
\int_{0}^{t}y^{\alpha -1}\phi _{\alpha }(y)\mathrm{d}y& =\int_{0}^{\infty
}u^{\alpha }G(\mathrm{d}u)\int_{0}^{t}y^{\alpha -1}e^{-uy}\mathrm{d}y \\
& =\int_{0}^{\infty }G(\mathrm{d}u)\int_{0}^{ut}z^{\alpha -1}e^{-z}\mathrm{d}%
z \\
& =\int_{0}^{\infty }\overline{G}(z/t)z^{\alpha -1}e^{-z}\mathrm{d}z \\
& =t^{\alpha }\int_{0}^{\infty }\overline{G}(u)u^{\alpha -1}e^{-ut}\mathrm{d}%
u.
\end{split}%
\end{equation*}%
A partial integration gives 
\begin{equation*}
1-\phi _{0}(t)=t\int_{0}^{\infty }\overline{G}(u)e^{-ut}\mathrm{d}u.
\end{equation*}%
So (\ref{eq:sconv}) reads 
\begin{equation}
t^{\alpha -1}\frac{\int_{0}^{\infty }\overline{G}(u)u^{\alpha -1}e^{-ut}%
\mathrm{d}u}{\int_{0}^{\infty }\overline{G}(u)e^{-ut}\mathrm{d}u}\rightarrow
\gamma \Gamma (\alpha ),\text{ as }t\searrow 0.  \label{eq:conv}
\end{equation}

From now on we shall assume that (\ref{gam}) holds with $0<\gamma \leq 1$.
Let us define the function for $t>0$ 
\begin{equation}
f(t)=\int_{0}^{\infty }\overline{G}(u)u^{\alpha -1}e^{-ut}\mathrm{d}u.
\label{eq:f-def}
\end{equation}%
Clearly, $f$ is monotone decreasing and since $\mathbb{E}Y=\infty $, $%
\lim_{t\rightarrow 0}f(t)=\infty $. Moreover, showing that $f$ is regularly
varying at zero implies that $\overline{G}$ is regularly varying at
infinity. We use the identity 
\begin{equation*}
u^{1-\alpha }e^{-ut}=\frac{1}{\Gamma (\alpha -1)}\int_{0}^{\infty }y^{\alpha
-2}e^{-(y+t)u}\mathrm{d}y,
\end{equation*}%
which holds for $u>0$ and $\alpha \in (1,2]$. (This is the \textit{%
Weyl-transform}, or \textit{Weyl-fractional integral }of the function $%
e^{-ut}$.) This identity combined with Fubini's theorem (everything is
nonnegative) gives 
\begin{equation*}
\begin{split}
\frac{1}{\Gamma (\alpha -1)}\int_{0}^{\infty }y^{\alpha -2}f(y+t)\mathrm{d}%
y& =\int_{0}^{\infty }\overline{G}(u)u^{\alpha -1}\mathrm{d}u\frac{1}{\Gamma
(\alpha -1)}\int_{0}^{\infty }y^{\alpha -2}e^{-(y+t)u}\mathrm{d}y \\
& =\int_{0}^{\infty }\overline{G}(u)e^{-ut}\mathrm{d}u.
\end{split}%
\end{equation*}%
So we can rewrite (\ref{eq:conv}) as 
\begin{equation}
\lim_{t \searrow 0}\frac{t^{\alpha -1}f(t)}{\int_{0}^{\infty }y^{\alpha
-2}f(t+y)\mathrm{d}y}=\frac{\gamma \Gamma (\alpha )}{\Gamma (\alpha -1)}%
=\gamma (\alpha -1).  \label{eq:conv-form-f1}
\end{equation}%
A change of variable gives 
\begin{equation*}
\int_{0}^{\infty }y^{\alpha -2}f(t+y)\mathrm{d}y=t^{\alpha
-1}\int_{1}^{\infty }(u-1)^{\alpha -2}f(ut)\mathrm{d}u,
\end{equation*}%
and so we have 
\begin{equation}
\lim_{t \searrow 0}\int_{1}^{\infty }(u-1)^{\alpha -2}\frac{f(ut)}{f(t)}%
\mathrm{d}u=[\gamma (\alpha -1)]^{-1}.  \label{eq:conv-form-f2}
\end{equation}%
We can rewrite $f$ as 
\begin{equation*}
f(t)=\int_{0}^{\infty }\overline{G}(u)u^{\alpha -1}e^{-ut}\mathrm{d}%
u=t^{-\alpha }\int_{0}^{\infty }\overline{G}(u/t)u^{\alpha -1}e^{-u}\mathrm{d%
}u,
\end{equation*}%
from which we see that the function 
\begin{equation*}
g(t)=\int_{0}^{\infty }\overline{G}(u/t)u^{\alpha -1}e^{-u}\mathrm{d}%
u=t^{\alpha }f(t)
\end{equation*}%
is bounded and nondecreasing. Substituting $g$ into (\ref{eq:conv-form-f2})
we obtain 
\begin{equation}
\lim_{t\rightarrow 0+}\int_{1}^{\infty }(u-1)^{\alpha -2}u^{-\alpha }\frac{%
g(ut)}{g(t)}\mathrm{d}u=[\gamma (\alpha -1)]^{-1}.  \label{eq:conv-form-g1}
\end{equation}

Write $g_{\infty }(x)=g(x^{-1})$, $x>0$. Then (\ref{eq:conv-form-g1}) has
the form 
\begin{equation}
\int_{1}^{\infty }(u-1)^{\alpha -2}u^{-\alpha }\frac{g_{\infty }(x/u)}{%
g_{\infty }(x)}\mathrm{d}u=\frac{k\overset{M}{\ast }g_{\infty }(x)}{%
g_{\infty }(x)}\rightarrow \lbrack \gamma (\alpha -1)]^{-1},\quad \text{as }%
x\rightarrow \infty ,  \label{eq:conv-form-g2}
\end{equation}%
where 
\begin{equation*}
k(u)=%
\begin{cases}
(u-1)^{\alpha -2}u^{-\alpha +1}, & u>1, \\ 
0, & 0<u\leq 1,%
\end{cases}%
\end{equation*}%
and 
\begin{equation*}
k\overset{M}{\ast }h(x)=\int_{0}^{\infty }h(x/u)k(u)/u\mathrm{d}u
\end{equation*}%
is the \textit{Mellin-convolution} of $h$ and $k$. Note that the \textit{%
Mellin-transform} of $k$, 
\begin{equation*}
\begin{split}
\widetilde{k}\left( z\right) & =\int_{1}^{\infty }\left( u-1\right) ^{\alpha
-2}u^{-\alpha -z}\mathrm{d}u =\int_{0}^{1}\left( 1-v\right) ^{\alpha -2}v^{z}%
\mathrm{d}v \\
& =\frac{\Gamma \left( \alpha -1\right) \Gamma \left( 1+z\right) }{\Gamma
\left( \alpha -z\right) }=\mathrm{Beta}(\alpha -1,1+z)
\end{split}%
\end{equation*}%
is convergent for $z>-1$. We apply a version of the Drasin-Shea theorem
(Theorem 5.2.3 on page 273 of \cite{BGT}). To do this we must verify the
following conditions:\medskip

\noindent1. $\widetilde {k}$ has a maximal convergent strip $a<\Re z<b$ such
that $a<0$ and $b>0,$ $\widetilde {k}\left( a+\right) =\infty$ and $%
\widetilde {k}\left( b-\right) =\infty$ if $b<\infty.$ Our $\widetilde {k}$
satisfies this condition with $a=-1$ and $b=\infty $. \medskip

\noindent2. Our function of interest is 
\begin{equation*}
g_{\infty}(x)=g(x^{-1})=\int_{0}^{\infty}\overline{G}(ux)u^{\alpha-1}e^{-u}%
\mathrm{d}u,\text{ }x>0,
\end{equation*}
is certainly positive and locally bounded. \medskip

\noindent3. Also our function $g_{\infty}$ is of bounded decrease, since for 
$\lambda>1$ 
\begin{equation*}
\frac{g_{\infty}(\lambda x)}{g_{\infty}(x)}=\lambda^{-\alpha}\frac{(\lambda
x)^{\alpha}g(1/(\lambda x))}{x^{\alpha}g(1/x)}=\lambda^{-\alpha}\frac{%
f(1/(\lambda x))}{f(1/x)}\geq\lambda^{-\alpha},
\end{equation*}
so its lower Matuszewska index is at least $-\alpha$.\medskip

Therefore by Theorem 5.2.3 of \cite{BGT}, whenever,%
\begin{equation}
\frac{k\overset{M}{\ast }g_{\infty }(x)}{g_{\infty }(x)}\rightarrow c,\quad 
\text{as }x\rightarrow \infty ,  \label{c}
\end{equation}%
then $\widetilde {k}(\rho )=c$ for some $\rho \in \left( -1,\infty \right) $%
. (In our case by (\ref{eq:conv-form-g2}), $c=[\gamma (\alpha -1)]^{-1}.$) \
Moreover, since $\widetilde{k}\left( z\right) $ is strictly decreasing on $%
\left( -1,\infty \right) $ and $\widetilde{k}\left( 0\right) =\frac{1}{%
\alpha -1},$ for any $0<\gamma \leq 1$ the solution $\rho $ to $\widetilde {k%
}(\rho )=[\gamma (\alpha -1)]^{-1}$ must lie in $(-1,0].$ Theorem 5.2.3 of 
\cite{BGT} also says that $g_{\infty }(x)$ is regularly varying at infinity
with index $0\geq \rho >-1$. \medskip

Next since $g_{\infty }(x)=g(x^{-1})=x^{-\alpha }f(x^{-1})\in \mathcal{RV}%
_{\infty }(\rho )$, where $\widetilde {k}(\rho )=c$, $g\in \mathcal{RV}%
_{0}(-\rho )$, which implies that $f\in \mathcal{RV}_{0}(-\rho -\alpha )$.
Recalling that%
\begin{equation*}
f(t)=\int_{0}^{\infty }\overline{G}(u)u^{\alpha -1}e^{-ut}\mathrm{d}u\text{,}
\end{equation*}%
the Karamata Tauberian theorem now gives that 
\begin{equation*}
\int_{0}^{x}\overline{G}(u)u^{\alpha -1}\mathrm{d}u\in \mathcal{RV}_{\infty
}(\alpha +\rho ).
\end{equation*}%
Thus by Lemma 2, $\overline{G}(u)\in \mathcal{RV}_{\infty
}(\rho).\medskip $

\noindent This says that $Y\in D(\beta )$, where $\rho =-\beta \in (-1,0]$
and $\beta $ is the unique solution of 
\begin{equation*}
\mathrm{Beta}(\alpha -1,\beta +1)=\frac{\Gamma (\alpha -1)\Gamma (1+\beta )}{%
\Gamma (\alpha -\beta )}=\frac{1}{\gamma (\alpha -1)}.
\end{equation*}

We now turn to the proof of the second part of Proposition 2. First consider
the the case $\beta =0$. Let $0 \leq D_{n}^{\left( n\right) }\leq \dots \leq
D_{n}^{\left( 1\right) }$ denote the order statistics of $Y_{1}/\left(
\sum_{i=1}^{n}Y_{i}\right), \dots, Y_{n}/\left(\sum_{i=1}^{n}Y_{i}\right)$.
We see that 
\begin{equation*}
\mathbb{E}\left( D_{n}^{\left( 1\right) }\right) ^{\alpha }\leq \mathbb{E}%
S_{n}\left( \alpha \right) =\sum_{i=1}^{n}\mathbb{E}\left( D_{n}^{\left(
i\right) }\right) ^{\alpha }\leq \mathbb{E}\left( D_{n}^{\left( 1\right)
}\right) ^{\alpha -1}\leq 1.
\end{equation*}%
Now $D_{n}^{\left( 1\right) }\rightarrow _{P}1$ if and only if $Y\in D(0)$.
(See Theorem 1 of Haeusler and Mason \cite{HM} and their references.) Thus
if $Y\in D(0) $ then (1.6) holds with $\gamma =1$.\medskip

Now assume that $Y\in D(\beta )$, $0<\beta <1$. In this case, there exists a
sequence of positive constants $\{a_{n}\}_{n\geq 1}$, such that $%
a_{n}^{-1}\sum_{i=1}^{n}Y_{i}\rightarrow _{d}U$, where $U$ is a $\beta $%
-stable random variable, with characteristic function 
\begin{equation*}
\mathbb{E}e^{\imath tU}=\exp \left\{ \beta \int_{0}^{\infty }(e^{\imath
tu}-1)u^{-\beta -1} \mathrm{d} u\right\} .
\end{equation*}%
Moreover, $Y^{\alpha }\in D(\beta /\alpha )$, and it is easy to check that $%
a_{n}^{-\alpha }\sum_{i=1}^{n}Y_{i}^{\alpha }\rightarrow _{d}V$, where $V$
is a $\beta /\alpha $-stable random variable, with cf 
\begin{equation*}
\mathbb{E}e^{\imath tV}=\exp \left\{ \frac{\beta }{\alpha }\int_{0}^{\infty
}(e^{\imath tu}-1)u^{-\beta /\alpha -1} \mathrm{d} u\right\} .
\end{equation*}%
Since 
\begin{equation*}
\lim_{n\rightarrow \infty }n\mathbb{P}\{Y>a_{n}u,Y^{\alpha }>a_{n}^{\alpha
}v\}=\lim_{n\rightarrow \infty }n\overline{G}(a_{n}(u\vee v^{1/\alpha
}))=u^{-\beta }\wedge v^{-\beta /\alpha }=:\Pi ((u,\infty ) \times (v,\infty )),
\end{equation*}%
for $u,v\geq 0$, $u+v>0$, using Corollary 15.16 of Kallenberg \cite{Kallen}
one can show that the joint convergence also holds, and the limiting
bivariate L\'{e}vy measure is $\Pi $. That is 
\begin{equation*}
\left( a_{n}^{-1}\sum_{i=1}^{n}Y_{i},a_{n}^{-\alpha
}\sum_{i=1}^{n}Y_{i}^{\alpha }\right) \rightarrow _{d}(U,V),
\end{equation*}%
where the limiting bivariate random vector has cf 
\begin{equation*}
\mathbb{E}e^{\imath (sU+tV)}=\exp \left\{ \int_{[0,\infty )^{2}}\left(
e^{\imath (su+tv)}-1\right) \Pi (\mathrm{d} u, \mathrm{d} v)\right\} =\exp \left\{ \beta
\int_{0}^{\infty }\left( e^{\imath (su+tu^{\alpha })}-1\right) u^{-\beta
-1}\mathrm{d} u\right\} .
\end{equation*}%
Since $\mathbb{P}\left\{ U>0\right\} =\mathbb{P}\left\{ V>0\right\} =1$, we
obtain 
\begin{equation*}
S_{n}\left( \alpha \right) \rightarrow _{d}\frac{V}{U^{\alpha }}.
\end{equation*}%
Thus since $\mathbb{E}S_{n}\left( \alpha \right) \leq 1$ for all $n\geq 1$%
\begin{equation*}
\mathbb{E}S_{n}\left( \alpha \right) \rightarrow \mathbb{E}\left( \frac{V}{%
U^{\alpha }}\right) .
\end{equation*}%
Clearly $\mathbb{P}\left\{ U < \infty \right\} = 1$, which implies that 
$0<\mathbb{E}\left( \frac{V}{U^{\alpha }}\right) \leq 1$,
and thus by the first part of Proposition 2  
\begin{equation*}
0<\gamma =\frac{\Gamma (\alpha -\beta )}{\Gamma (\alpha )\Gamma (1+\beta )}%
<1.
\end{equation*}%
\hfill \mbox{$\Box$}

\noindent \textbf{Lemma 2} \textit{Suppose that for some $\alpha \geq 1$, $%
\rho >-1$ and slowly varying function $L$ at infinity 
\begin{equation*}
U\left( x\right) :=\int_{0}^{x}\overline{G}(u)u^{\alpha -1}\mathrm{d}%
u=L\left( x\right) x^{\alpha +\rho },\text{ }x>0,
\end{equation*}%
then 
\begin{equation*}
\overline{G}(u)\sim \left( \alpha +\rho \right) L\left( u\right) u^{\rho },%
\text{ as }u\rightarrow \infty .
\end{equation*}%
}

\noindent \textit{Proof} We shall follow closely the proof the lemma on page
446 of Feller \cite{Feller}. Choose any $0<a<b<\infty $. We see that 
\begin{equation*}
\frac{U\left( tb\right) -U\left( ta\right) }{U\left( t\right) }=\int_{a}^{b}%
\frac{\overline{G}(ut)\left( ut\right) ^{\alpha -1}t}{U\left( t\right) }%
\mathrm{d}u
\end{equation*}%
\begin{equation*}
=\int_{a}^{b}\frac{\overline{G}(ut)\left( ut\right) ^{\alpha -1}t}{L\left(
t\right) t^{\alpha +\rho }}\mathrm{d}u=\int_{a}^{b}\frac{\overline{G}%
(ut)u^{\alpha -1}\mathrm{d}u}{L\left( t\right) t^{\rho }}.
\end{equation*}%
Since $\overline{G}$ is nonincreasing and $\overline{G}\left( ut\right)
/\left( L\left( t\right) t^{\rho }\right) $ is necessarily bounded for each $%
u>0$ as $t\rightarrow \infty $, just as in Feller one can apply the
Helly-Bray theorem to find a positive sequence $t_{k}\rightarrow \infty $
such that for a measurable function $\psi $ on $\left[ 0,\infty \right) $, $%
\overline{G}(ut_{k})/L\left( t_{k}\right) t_{k}^{\rho }\rightarrow $ $\psi
\left( u\right) $, for all continuity points $u$ of $\psi $. This implies
that for all $0<a<b<\infty $%
\begin{equation*}
\frac{U\left( t_{k}b\right) -U\left( t_{k}a\right) }{U\left( t_{k}\right) }%
\rightarrow b^{\alpha +\rho }-a^{\alpha +\rho }=\int_{a}^{b}\psi \left(
u\right) u^{\alpha -1}\mathrm{d}u.
\end{equation*}%
This forces $\psi \left( u\right) u^{\alpha -1}=\left( \alpha +\rho \right)
u^{\alpha +\rho -1}$, and since $\psi $ is independent of any particular
positive sequence $t_{k}\rightarrow \infty $ defining it, 
\begin{equation*}
\overline{G}(ut)/\left( L\left( ut\right) \left( ut\right) ^{\rho }\right)
\rightarrow \alpha +\rho \text{, as } t \rightarrow \infty .
\end{equation*}%
\hfill\mbox{$\Box$}

\medskip

\noindent \textbf{Acknowledgement.}
Kevei's research was funded by a postdoctoral fellowship
of the Alexander von Humboldt Foundation.

\end{document}